\documentclass[12pt,dvipdfmx]{article}
\usepackage{amsmath}
\usepackage{tikz}
\usetikzlibrary{intersections,calc,arrows}
\usepackage{amssymb}
\date{}
\makeatletter
\newcommand{\figcaption}[1]{\def\@captype{figure}\caption{#1}}
\newcommand{\tblcaption}[1]{\def\@captype{table}\caption{#1}}

\@addtoreset{equation}{section}
\makeatother
\newcommand{\qed}{\hbox{\rule[-2pt]{3pt}{6pt}}}

\begin{document}
\title {\bf Bifurcation diagrams of one-dimensional Kirchhoff type equations}

\author{Tetsutaro Shibata 
\\
Laboratory of Mathematics, School of Engineering
\\
Graduate School of Advanced Science and Engineering
\\
Hiroshima University, 
Higashi-Hiroshima, 739-8527, Japan
}

\maketitle
\footnote[0]{E-mail: tshibata@hiroshima-u.ac.jp}
\footnote[0]{This work was supported by JSPS KAKENHI Grant Number JP21K03310.}

\begin{abstract}
We study the one-dimensional Kirchhoff type equation
$$
-(b + a\Vert u'\Vert^{2}) u''(x) = \lambda u(x)^p, x \in I:= (-1,1), \enskip u(x) > 0, 
\enskip x\in I, \enskip u(\pm 1) = 0,
$$
where $\Vert u'\Vert = \left(\int_I u'(x)^2 dx\right)^{1/2}$, $a > 0, b > 0, p> 0$ are given constants and $\lambda > 0$ is a bifurcation parameter.  
We establish the exact solution $u_\lambda(x)$ and complete shape of the bifurcation curves $\lambda = \lambda(\xi)$, where 
$\xi:= \Vert u_\lambda\Vert_\infty$. 
We also study the nonlinear eigenvalue problem
$$
-\Vert u'\Vert^{p-1} u''(x) = \mu u(x)^p, x \in I, \enskip u(x) > 0,  x\in I, \enskip u(\pm 1) = 0,
$$
where $p > 1$ is a given constant and $\mu > 0$ is an eigenvalue parameter. 
We establish the first eigenvalue and eigenfunction of this problem by using a simple time map method. 

\end{abstract}


\noindent
{Keywords:}

\section{Introduction} 		      

We study the structure of the global solution curves for nonlocal elliptic problem

\begin{equation}
\left\{ \,
    \begin{aligned}
    -(b + a\Vert u'\Vert^{2}) u''(x) &= \lambda u(x)^p, \enskip x \in I:= (-1,1),  \\
    u(x) &> 0, \enskip x\in I,  \\
    u(\pm 1) &= 0. 
    \end{aligned}
\right.
\end{equation}
where $\Vert u'\Vert = \left(\int_I u'(x)^2 dx\right)^{1/2}$, $a > 0, b > 0, p > 0$ are given constants. 
 If $p = 3$, then the problem (1.1) is known as the one-dimensional elliptic Kirchhoff type equation. 

Nonlinear elliptic bifurcation problems have been studied intensively by many authors. 
We refer to [5] and the references therein. Besides,  
nonlocal elliptic problems have been also studied by many investigators, since they are 
derived from several interesting physical and engineering phenomena. We refer to [1] and the references therein.  In the field of nonlocal elliptic problems, 
there are several studies which deal with bifurcation problems. 
We refer to [1,3,4,8,9] and the references therein. 
As for as the author knows, however, there are few results which clarify the 
precise structures of bifurcation diagrams for nonlocal problems.  Our problem (1.1) was proposed in [11] as the elliptic eigenvalue problems 
in bounded smooth domain $\Omega \subset \mathbb{R}^n$ ($n=1,2,3$)
\begin{equation}
\left\{ \,
    \begin{aligned}
-\left(b + a\int_\Omega \vert \nabla u\vert^2 dx\right) \Delta u &= \lambda u^p \quad  \mbox{in} \enskip  \Omega, 
\\
u &> 0  \quad \mbox{in} \enskip  \Omega,
\\
u &= 0 \quad \mbox{on} \enskip  \partial\Omega,
\end{aligned}
\right.
\end{equation}
and using the degree argument and variational method, the following results have been obtained.

\vspace{0.2cm}

\noindent
{\bf Theorem 1.1 ([11]).} {Consider (1.2). \it Let $2^*= \infty$ for $n = 1,2$ and $2^* = 6$ for $n = 3$.  

\noindent
{(i)} Assume that $0 < p < 1$. Then (1.2) has a branch of positive solutions bifurcating from 
zero at $\lambda = 0$. 

\noindent
{(ii)} Assume that $3 < p < 2^*-1$. Then (1.2) has a branch of positive solutions bifurcating from infinity at $\lambda = 0$. 
}

\vspace{0.2cm}
 
It should be mentioned that, in many cases, $\lambda$ is parameterized by the maximum 
norm $\xi:=\Vert u_\lambda \Vert_\infty$ 
as $\lambda = \lambda(\xi)$, where $u_\lambda$ is a solution associated with $\lambda$.  
However, $\lambda$ is not expressed explicitly by using $\xi$ in [11]. 
Further, in [11], the case $1 \le p \le 3$ was 
not considered. 
\noindent
 In this paper, we focus on the case $n = 1$ and establish the exact representation of 
 $\lambda(\xi)$ for $p > 0$ by using simple time map method. We put 
\begin{eqnarray}
A_p &:=& \int_0^1 \frac{1}{\sqrt{1-s^{p+1}}}ds,
\\
B_p &:=& \int_0^1 \sqrt{1-s^{p+1}}ds,
\\
C_p &:=& \int_0^1 \frac{s^{p+1}}{\sqrt{1-s^{p+1}}}ds.
\end{eqnarray}
Let $W_p(x)$ ($0 < p < 1$, $p > 1$) be a unique solution of 
\begin{equation}
\left\{ \,
    \begin{aligned}
-W''(x) &= W(x)^p, \quad x \in I,
\\
W(x) &> 0, \quad x \in I,
\\
W(-1) &= W(1) = 0.
\end{aligned}
\right.
\end{equation}
The unique existence of $W_p(x)$ is known by [2, 6]. 

Now we state the first results.

\vspace{0.2cm}

\noindent
{\bf Theorem 1.2.} {\it Consider (1.1). Let 
$a > 0, b > 0$ and $p > 0$ be given constants. 
Assume that $0 < p < 1$ or $3 < p < \infty$. 
Then for any $\lambda > 0$, there exists a unique solution pair $(\lambda, u_\lambda) \in 
\mathbb{R}_+ \times C^2(\bar{I})$. Put $\xi:= \Vert u_\lambda\Vert_\infty$. Then 
$\lambda$ is parameterized by $\xi$ and the following 
formula holds for $\xi > 0$. 
\begin{eqnarray}
\lambda(\xi) = \frac{p+1}{2}A_p^2\left(2A_pB_p a\xi^{3-p} + b\xi^{1-p}\right).
\end{eqnarray}
Furthermore, $u_\lambda(x)$ is given by
\begin{eqnarray}
u_\lambda(x) = \left(\frac{\lambda}{\sqrt{b + 2A_pB_pa\xi(\lambda)^2}}\right)^{1/(1-p)}
W_p(x),
\end{eqnarray}
where $\xi(\lambda)$ is the inverse function of $\lambda(\xi)$. 
}

\vspace{0.2cm}

\noindent
We note that, if $0 < p < 1$ (resp. $3 < p < \infty$), then we see from (1.7) that 
$\lambda(\xi)$ is strictly increasing (resp. decreasing) for $\xi > 0$. Therefore, $\xi(\lambda)$ 
exists.

\vspace{0.2cm}

Theorem 1.2 improves the results in Theorem 1.1 and gives the explicit formula for 
bifurcation curves of the problem (1.1) for the case $n = 1$. 

Next, we consider the case where $1 < p < 3$. We apply [1, Theorem 2] to (1.1) and obtain the following Theorem 1.3. 

\vspace{0.2cm}

\noindent
{\bf Theorem 1.3.} {\it Consider (1.1). Let 
$a > 0, b > 0$ and $1 < p < 3$ be given constants. 
For $\lambda > 0$, put 
\begin{eqnarray}
L(\lambda):=\frac{p-1}{2}\lambda^{2/(p-1)}\Vert W_p'\Vert^{-2}
\left(\frac{2b}{3-p}\right)^{(p-3)/(p-1)}. 
\end{eqnarray}
{(a)} If $L(\lambda) > a$, then (1.1) has exactly two solutions. 

\noindent
{(b)} If $L(\lambda) = a$, then (1.1) has exactly one solution. 

\noindent
{(c)} If $L(\lambda) < a$, then (1.1) has no solutions. 
}

\vspace{0.2cm}

\noindent
Unfortunately, it is rather difficult to obtain the clear shape of $\lambda(\xi)$ for 
general 
$1 < p < 3$. The reason will be explained in Section 4. 
Therefore, we concentrate on the special case $p = 2$ and establish the formulas for $\xi = \xi(\lambda)$. We put 
\begin{eqnarray}
Q_1&:=& \frac{3}{2}A_2^2,
\\
Q_2&:=& \frac{3}{\sqrt{2}}A_2^{5/2}B_2^{1/2}.
\end{eqnarray}

\vspace{0.2cm}

\noindent
{\bf Theorem 1.4.} {\it Consider (1.1). Let $p = 2$ and 
$a > 0, b > 0$ be given constants. 

\noindent
{(i)} Let 
\begin{eqnarray}
\lambda = 2\sqrt{ab}Q_2. 
\end{eqnarray}
Then (1.1) has a unique solution $u_\lambda(x) = 2b\lambda^{-1}W_2(x)$. Moreover,
\begin{eqnarray}
\xi = 2bQ_1\lambda^{-1} = \sqrt{\frac{b}{a}}\frac{Q_1}{Q_2}.
\end{eqnarray}

\noindent
{(ii)} Let $\lambda > 2\sqrt{ab}Q_2$. Then there exist exactly two solutions
$u_{\lambda,1}(x)$ and $u_{\lambda,2}(x)$. Moreover, $\xi$ are parameterized 
by $\lambda$ such as 
$\xi_1= \xi_1(\lambda)$, $\xi_2 = \xi_2(\lambda)$, and are represented as follows. 
\begin{eqnarray}
\xi_1(\lambda)& =& \frac{\lambda Q_2^{-1} + \sqrt{\lambda^{2}Q_2^{-2}-4ab}}{2a}Q_2^{-1}Q_1,
\\
\xi_2(\lambda)& =& \frac{\lambda Q_2^{-1} - \sqrt{\lambda^{2}Q_2^{-2}-4ab}}{2a}Q_2^{-1}Q_1.
\end{eqnarray}

\noindent
{(iii)} Let $\lambda < 2\sqrt{ab}Q_2$. Then there are no solutions of (1.1). 

}

\vspace{0.2cm}

Before considering the case $p = 1$ and $p = 3$ for (1.1),  we study the 
following nonlinear eigenvalue problems 

\begin{equation}
\left\{ \,
    \begin{aligned}
-\Vert u'\Vert^{p-1} u''(x) &= \mu u(x)^p, \enskip x \in I, 
\\
u(x) &> 0, \enskip x\in I, 
\\
u(\pm 1) &= 0,
\end{aligned}
\right.
\end{equation}
where 
$p > 1$ is a given constant,  and $\mu > 0$ is 
an eigenvalue parameter. In the case 
$p = 3$ in (1.15), it is also known as Kirchhoff type eigenvalue problem. 
We refer to [3], which treated the case 
$p = 3$ precisely, and the refernces therein. 
Unfortunately, however, the exact number of $\mu_1$ was not mentioned even in [3]. 
Here, we obtain the explicit first eigenvalue $\mu_1$ of the 
problem (1.15) by using 
simple time map method, where $\mu_1$ is defined by 
\begin{eqnarray}
\mu_1 := \inf\left\{\Vert u'\Vert^{p+1}: u \in H_0^1(I), \int_I u(x)^{p+1}dx = 1\right\}.
\end{eqnarray}
We also obtain the explicit form of the first eigenfunction 
$\varphi_1(x)$ associated with $\mu_1$. The existence of $\mu_1$ and $\varphi_1$ can be proved 
easily by the standard direct method used in the proof of the case $p = 1$. Namely, we choose 
the minimizing sequence $\{u_n\}_{n=1}^\infty \subset H_0^1(I)$ and use the compact embedding 
$H_0^1(I) \subset C(\bar{I})$ to obtain the minimizer $\varphi_1(x)$ and $\mu_1$.

\vspace{0.2cm}

\noindent
{\bf Theorem 1.5.} {\it Consider (1.15). 
Then $\mu_1$ is given by 
\begin{eqnarray}
\mu_1 = 2^{(p-3)/2}(p+1)A_p^{(p+3)/2}B_p^{(p-1)/2}.
\end{eqnarray}
Furthermore, 
\begin{eqnarray}
\varphi_1(x) =  2^{2/(p^2-1)}(p+1)^{-1/(p-1)}A_p^{-(p+3)/(p^2-1)}C_p^{-1/(p+1)}W_p(x). 
\end{eqnarray}
}

\vspace{0.2cm}

Now we consider the case $p = 1$ and $p = 3$ for (1.1). For $p = 3$, we know from [10, Theorem 1.2] that 
the following result holds true.  

\vspace{0.2cm}

\noindent
{\bf Theorem 1.6 ([10, Theorem 1.2]).} {\it Assume that $p = 3$. If $a > 0, b > 0$ and 
$\lambda > a\mu_1$, then (1.1) has at least one positive solution. }

\vspace{0.2cm}

We improve Theorem 1.6  for the case $n = 1$. 

\vspace{0.2cm}

\noindent
{\bf Theorem 1.7.} {\it (i) Assume that $p = 3$. Let $\lambda > a\mu_1= 4aA_3^3B_3$. Then 
(1.1) has a unique solution 
\begin{eqnarray}
u_\lambda(x) =  \sqrt{\frac{b\Vert W_3'\Vert^2}{\lambda - a\Vert W_3'\Vert^2}}
\Vert W_3'\Vert^{-1}W_3(x),
\end{eqnarray} 
where $\Vert W_3'\Vert = 2A_3^{3/2}B_3^{1/2}$ and  
\begin{eqnarray}
\lambda(\xi) = 4aA_3^3B_3 + 2A_3^2b\xi^{-2} \quad (\xi > 0).
\end{eqnarray}
(ii) Assume that $p = 1$. Let $\lambda > \frac{\pi^2}{4}b$. Then the solution $u_\lambda(x)$ 
of (1.1) is given as 
follows.
\begin{eqnarray}
u_\lambda(x) = \frac{4}{\pi^2}\sqrt{\frac{\lambda-(\pi^2b)/4}{a}}\cos\frac{\pi}{2}x.
\end{eqnarray}
Furthermore, 
\begin{eqnarray}
\lambda(\xi) = \frac{\pi^4}{16}a\xi^2 + \frac{\pi^2}{4}b \quad (\xi > 0).
\end{eqnarray}
}

\vspace{0.2cm}

The remainder of this paper is as follows. In Section 2, we prove Theorem 1.5 to introduce 
the time map method (cf. [7]). By using this argument, Theorems 1.2 and 1.4 will be proved in Sections 2 an 3, respectively.  Finally, we will prove Theorem 1.7 in Section 5.

\section{Proof of Theorem 1.5}

We first prove (1.17). We have 
\begin{eqnarray}
\varphi_1(x) &=& \varphi_1(1-x), 
\\
\zeta&:=& \max_{x \in I}\varphi_1(x) = \varphi_1(0),
\\
\varphi_1'(x) &>& 0, \enskip -1 \le x < 0.
\end{eqnarray}
We introduce the time map argument. From (1.15), we have 
\begin{eqnarray}
\{\Vert \varphi_1'\Vert^{p-1} \varphi_1''(x) + \mu_1 \varphi_1(x)^{p}\}\varphi_1'(x) = 0.
\end{eqnarray}
This implies that 
\begin{eqnarray}
\frac{d}{dx}\left\{\frac12\Vert \varphi_1'\Vert^{p-1} \varphi_1'(x)^2 + \frac{1}{p+1}\mu_1 
\varphi_1(x)^{p+1}\right\} = 0.
\end{eqnarray}
By puting $x = 0$ and (2.2), for $x \in I$, we have 
\begin{eqnarray}
\frac12\Vert \varphi_1'\Vert^{p-1} \varphi_1'(x)^2 + \frac{1}{p+1}\mu_1 \varphi_1(x)^{p+1} &=& \mbox{constant} 
= \frac{1}{p+1} \mu_1 \zeta^{p+1}.
\end{eqnarray}
By this and (2.3), for $-1 \le x < 0$, we have 
\begin{eqnarray}
\varphi_1'(x) = \sqrt{\frac{2}{p+1}}\Vert \varphi_1'\Vert^{-(p-1)/2}\sqrt{\mu_1(\zeta^{p+1} 
- \varphi_1(x)^{p+1})}.
\end{eqnarray}
This along with (1.3) implies that 
\begin{eqnarray}
\sqrt{\mu_1} &=& \sqrt{\frac{p+1}{2}}\Vert \varphi_1'\Vert^{(p-1)/2}\int_{-1}^0 
\frac{\varphi_1'(x)}{\sqrt{\zeta^{p+1}-\varphi_1(x)^{p+1}}}dx
\\
&=& \sqrt{\frac{p+1}{2}}\Vert \varphi_1'\Vert^{(p-1)/2}
\int_0^\zeta \frac{1}{\sqrt{\zeta^{p+1}-\theta^{p+1}}}d\theta
\nonumber
\\
&=& \sqrt{\frac{p+1}{2}}\Vert \varphi_1'\Vert^{(p-1)/2}\zeta^{(1-p)/2}
\int_0^1 \frac{1}{\sqrt{1-s^{p+1}}}ds 
\nonumber
\\
&=& \sqrt{\frac{p+1}{2}}\Vert \varphi_1'\Vert^{(p-1)/2}\zeta^{(1-p)/2}A_p.
\nonumber
\end{eqnarray}
By using (1.5), (1.16), (2.1) and (2.7), we have 
\begin{eqnarray}
\frac12 &=& \int_{-1}^0 \varphi_1(x)^{p+1}dx 
\\
&=& \sqrt{\frac{p+1}{2}}\Vert \varphi_1'\Vert^{(p-1)/2}
\frac{1}{\sqrt{\mu_1}}\int_{-1}^0 \frac{\varphi_1(x)^{p+1}\varphi_1'(x)}{\sqrt{\zeta^{p+1}-\varphi_1(x)^{p+1}}}dx 
\nonumber
\\
&=& \sqrt{\frac{p+1}{2}}\Vert \varphi_1'\Vert^{(p-1)/2}
\frac{1}{\sqrt{\mu_1}}\int_0^\zeta \frac{\theta^{p+1}}{\sqrt{\zeta^{p+1}-\theta^{p+1}}}d\theta
\nonumber
\\
&=& \sqrt{\frac{p+1}{2}}\Vert \varphi_1'\Vert^{(p-1)/2}
\frac{1}{\sqrt{\mu_1}}\zeta^{(p+3)/2} C_p.
\nonumber
\end{eqnarray}
By (2.8) and (2.9), we obtain 
\begin{eqnarray}
\zeta = \left(\frac{A_p}{2C_p}\right)^{1/(p+1)}.
\end{eqnarray}
By (1.4), (2.7) and (2.8), we have 
\begin{eqnarray}
\Vert \varphi_1'\Vert^2 &=& 2\int_{-1}^0 \varphi_1'(x)^2 dx 
\\
&=& 2\sqrt{\frac{2}{p+1}}\Vert \varphi_1'\Vert^{-(p-1)/2}\sqrt{\mu_1}
\int_{-1}^0\sqrt{\zeta^{p+1}-\varphi_1(x)^{p+1}}\varphi_1'(x)dx
\nonumber
\\
&=& 2A_pB_p\zeta^2.
\nonumber
\end{eqnarray}
This along with (2.10) implies that 
\begin{eqnarray}
\Vert \varphi_1'\Vert = \sqrt{2A_pB_p}\left(\frac{A_p}{2C_p}\right)^{1/(p+1)}.
\end{eqnarray}
By this, (2.8) and (2.10), we obtain
\begin{eqnarray}
\mu_1 = \frac{p+1}{2}2^{(p-1)/2}A_p^{(p+3)/2}B_p^{(p-1)/2}.
\end{eqnarray}
Thus we obtain (1.17). We next prove (1.18). By (1.15), (1.17) and (2.12), we obtain 
\begin{eqnarray}
-\varphi_1''(x) = 2^{-2/(p+1)}(p+1)A_p^{(p+3)/(p+1)}C_p^{(p-1)/(p+1)}\varphi_1(x)^{p}.
\end{eqnarray}
We put 
\begin{eqnarray}
\nu:= 2^{2/(p^2-1)}(p+1)^{-1/(p-1)}A_p^{-(p+3)/(p^2-1)}C_p^{-1/(p+1)}, \enskip \varphi_1(x) = \nu W_p(x).
\end{eqnarray}
Then we see that $\varphi_1(x)$ satisfies (1.15) with $\mu_1$. 
This implies (1.18). Thus the proof of Theorem 1.5 is complete. \qed

\section{Proof of Theorem 1.2} Let $p > 0$ ($p \not= 1$). 
We consider the following equation of $t \in \mathbb{R}$. 
\begin{eqnarray}
at + b = \lambda\Vert W_p'\Vert^{1-p} t^{(p-1)/2}.
\end{eqnarray}
Assume that there exists a solution
$t_\lambda > 0$ of (3.1). Then we see from [1, Theorem 2] that there exists a solution pair 
$(u_\lambda,\lambda)$ of (1.1) corresponding to $t_\lambda$, and $u_\lambda(x)$ is given by 
\begin{eqnarray}
u_\lambda(x) = t_\lambda^{1/2}\Vert W_p'\Vert^{-1}W_p(x).
\end{eqnarray}
Indeed, let $w_\lambda$ be a unique solution of 
\begin{equation}
\left\{ \,
    \begin{aligned}
-w''(x) &= \lambda w(x)^p, \quad x \in I,
\\
w(x) &> 0, \quad x \in I,
\\
w(-1) &= w(1) = 0.
\end{aligned}
\right.
\end{equation}
Then we see that $w_\lambda = \lambda^{1/(1-p)}W_p$. We put $\gamma = 
t_\lambda^{1/2}\Vert w_\lambda'\Vert^{-1}$ and $u_\lambda:= \gamma w_\lambda$, 
Then we have 
\begin{eqnarray}
b + a\Vert \gamma w_\lambda'\Vert^2 = b + at_\lambda = t_\lambda^{(p-1)/2}
\Vert w_\lambda'\Vert^{1-p} = \gamma^{p-1}.
\end{eqnarray}
Then we have 
\begin{eqnarray}
-(b + a\Vert u_\lambda'\Vert^2)u_\lambda'' = -\gamma^pw_\lambda'' 
= \lambda \gamma^p w_\lambda^p = \lambda u_\lambda.
\end{eqnarray}
We note that
\begin{eqnarray}
u_\lambda = \gamma w_\lambda = t_\lambda^{1/2}\lambda^{1/(p-1)}\Vert W_p'\Vert^{-1}
\lambda^{1/(1-p)}W_p =  t_\lambda^{1/2}\Vert W_p'\Vert^{-1}W_p.
\end{eqnarray}
Therefore, if (3.1) has $k$ positive solutions $t_{\lambda, 1,}, t_{\lambda,2}, \cdots,t_{\lambda,k}$, 
then (1.1) has also $k$ solutions correspoinging to $t_{\lambda,j}$ ($j = 1,2,\cdots,k$). 
On the contrary, assume that $u_\lambda$ is a solution of (1.1). We put 
$t_\lambda:= \Vert u_\lambda'\Vert^2$. Then we see that 
$$
u_\lambda(x) = \left(\frac{b + at_\lambda}{\lambda}\right)^{1/(p-1)}W_p(x), \quad
t_\lambda = \Vert u_\lambda'\Vert^2 =  \left(\frac{b + at_\lambda}{\lambda}\right)^{2/(p-1)}
\Vert W_p'\Vert^2.
$$
This implies that $t_\lambda$ satisfies (3.1). 
Therefore, the solutions of (1.1) correspond to those of  (3.1). 
For $t > 0$, we put 
\begin{eqnarray}
g(t):= at + b - Rt^{(p-1)/2},
\end{eqnarray}
where $R:= \lambda\Vert W_p'\Vert^{1-p}$. Now we look for the solutions of $g(t) = 0$. We have 
\begin{eqnarray}
g'(t) = a - \frac{p-1}{2}Rt^{(p-3)/2}. 
\end{eqnarray}
If $p > 3$, then $g(t)$ attains its maximum at 
\begin{eqnarray}
t_0:= \left(\frac{2a}{(p-1)R}\right)^{2/(p-3)} = 
\left(\frac{2a}{(p-1)\lambda\Vert W_p'\Vert^{1-p}}\right)^{2/(p-3)}. 
\end{eqnarray}
Since $g(0) = b$ and $g(t)$ strictly increases in $(0, t_0)$ and attains its maximum $g(t_0) > 0$ 
at $t_0$ 
and strictly decreases in $(t_0, \infty)$ and $g(\infty) = -\infty$. So it is clear that there 
exists a unique $t_\lambda > t_0$ which satisfies $g(t_\lambda) = 0$. 
If $0 < p < 1$, then by (3.7) and (3.8), we have $g'(t) > 0$ for $t > 0$ and $g(0) = -\infty$ and $g(\infty) = \infty$. 
So it is clear that there 
exists a unique $t_\lambda > 0$ which satisfies $g(t_\lambda) = 0$. 
Certainly, if $t_\lambda$ is represented explicitly by $\lambda$, then it is natural to find 
the relationship between $\lambda$ and $\xi = \Vert u_\lambda\Vert_\infty$, and it seems 
possible to 
obtain the bifurcation curve $\lambda(\xi)$. Unfortunately, however, it is difficult to 
obtain $t_\lambda$ explicitly. To overcome this difficulty, we apply time map method to obtain the 
bifurcation 
curves $\lambda(\xi)$ of (1.1). 

Let $p > 3$ or $0 < p < 1$ be fixed. We apply time map method to (1.1). Let an arbitrary $\lambda > 0$ be fixed and $u_\lambda(x)$ be a unique solution 
of (1.1). 
Then we have 
\begin{eqnarray}
u_\lambda(x) &=& u_\lambda(-x), \quad x \in [-1, 0],
\\
\xi &:=& \Vert u_\lambda\Vert_\infty = \max_{-1 \le x \le 1}u_\lambda(x) = u_\lambda(0),
\\
u_\lambda'(x) &>& 0, \quad x \in [-1,0).
\end{eqnarray}
By (1.1), we have 
\begin{eqnarray}
\{(b + a\Vert u'\Vert^2)u''(x) + \lambda u(x)^p\}u'(x) = 0.
\end{eqnarray}
This implies that for $-1 \le x \le 0$,
\begin{eqnarray}
(b + a\Vert u'\Vert^2)u'(x)^2 + \frac{2}{p+1}\lambda u_\lambda(x)^{p+1} 
= \frac{2}{p+1}\lambda \xi^{p+1}.
\end{eqnarray}
By this, for $-1 \le x \le 0$, we obtain 
\begin{eqnarray}
u'_\lambda(x) &=& \sqrt{\frac{2}{p+1}}\frac{1}{\sqrt{b + a\Vert u_\lambda'\Vert^2}}
\sqrt{\lambda}\sqrt{\xi^{p+1}-u_\lambda(x)^{p+1}}
\end{eqnarray}
 By this and (1.3), we obtain 
 \begin{eqnarray}
 \sqrt{\lambda} &=& \sqrt{\frac{p+1}{2}}\sqrt{b + a\Vert u_\lambda'\Vert^2}\int_{-1}^0
 \frac{u_\lambda'(x)}{\sqrt{\xi^{p+1}-u_\lambda(x)^{p+1}}}dx
 \\
 &=& \sqrt{\frac{p+1}{2}}\sqrt{b + a\Vert u_\lambda'\Vert^2}\int_0^{\xi}
 \frac{1}{\sqrt{\xi^{p+1}-\theta^{p+1}}}d\theta
 \nonumber
 \\
 &=& \sqrt{\frac{p+1}{2}}\sqrt{b + a\Vert u_\lambda'\Vert^2}\xi^{(1-p)/2}
 \int_0^1 \frac{1}{\sqrt{1-s^{p+1}}}ds
 \nonumber
 \\
 &=&\sqrt{\frac{p+1}{2}}\sqrt{b + a\Vert u_\lambda'\Vert^2}\xi^{(1-p)/2}A_p.
 \nonumber
 \end{eqnarray}
 By (1.4) and (3.15), we have 
 \begin{eqnarray}
 \Vert u_\lambda'\Vert^2 &=& 2\int_{-1}^0 u_\lambda'(x)u_\lambda'(x)dx 
 \\
 &=& 2\int_{-1}^0 \sqrt{\frac{2}{p+1}}\frac{1}{\sqrt{b + a\Vert u'\Vert^2}}
 \sqrt{\lambda}\sqrt{\xi^{p+1}-u_\lambda(x)^{p+1}}
 u_\lambda'(x)dx 
 \nonumber
 \\
 &=& 2\sqrt{\frac{2}{p+1}}\frac{1}{\sqrt{b + a\Vert u_\lambda'\Vert^2}}
 \sqrt{\lambda}\int_0^\xi \sqrt{\xi^{p+1}-\theta^{p+1}}dx
 \nonumber
 \\
 &=& 2\sqrt{\frac{2}{p+1}}\frac{\sqrt{\lambda}}{\sqrt{b + a\Vert u'\Vert^2}}
 \xi^{(p+3)/2}\int_0^1 \sqrt{1-s^{p+1}}ds 
 \nonumber
 \\
 &=& 
2\sqrt{\frac{2}{p+1}}\frac{\sqrt{\lambda}}{\sqrt{b + a\Vert u_\lambda'\Vert^2}}
 \xi^{(p+3)/2}B_p.
 \nonumber
 \end{eqnarray}
 Substitute (3.16) into (3.17). Then we have 
 \begin{eqnarray}
 \Vert u_\lambda'\Vert^2 &=& 2A_pB_p\xi^2.
 \end{eqnarray}
 By this and (3.16), we have 
 \begin{eqnarray}
 \sqrt{\lambda} &=& \sqrt{\frac{p+1}{2}}\sqrt{b + 2A_pB_pa\xi^2}\xi^{(1-p)/2}A_p.
 \end{eqnarray}
 Namely, 
 \begin{eqnarray}
 \lambda= \frac{p+1}{2}A_p^2(2A_pB_pa\xi^{3-p} + b\xi^{1-p}).
 \end{eqnarray}
 This implies (1.7). We next prove (1.8). By (1.1), (1.7) and (3.18), we have 
 \begin{eqnarray}
 -u_\lambda''(x) &=& \frac{\lambda}{\sqrt{b + a\Vert u_\lambda'\Vert^2}}u_\lambda(x)^p
 \\
 &=& \frac{\lambda}{\sqrt{b + 2A_pB_pa\xi(\lambda)^2}}u_\lambda(x)^p.
 \nonumber
 \end{eqnarray}
 This implies that 
 \begin{eqnarray}
 u_\lambda(x) = \left(\frac{\lambda}{\sqrt{b + 2A_pB_pa\xi(\lambda)^2}}\right)^{1/(1-p)}
 W_p(x).
 \end{eqnarray}
 This implies (1.8). Thus the proof of Theorem 2 is complete. \qed 

\section{Proof of Theorems 1.3 and 1.4} 

We begin with the proof of Theorem 1.3. 

\vspace{0.2cm}

\noindent
{\bf Proof of Theorem 1.3.} Let $1 < p < 3$. We fix $b > 0$. Assume that $y(t) = at + b$ 
is tangent line of 
$g(t) = \lambda\Vert W_p'\Vert^{1-p}t^{(p-1)/2}$.  Let $t = t_0$ be a point of tangency. Then we have 
\begin{eqnarray}
a = \frac{p-1}{2}\lambda \Vert W_p'\Vert^{1-p}t_0^{(p-3)/2}.
\end{eqnarray}
Since $at_0 + b = g(t_0)$, by (4.1), we have 
\begin{eqnarray}
at_0 + b = \frac{p-1}{2}\lambda \Vert W_p'\Vert^{1-p}t_0^{(p-1)/2} + b = \lambda
\Vert W_p'\Vert^{1-p}t_0^{(p-1)/2}.
\end{eqnarray}
By this, we have 
\begin{eqnarray}
b = \frac{3-p}{2}\lambda
\Vert W_p'\Vert^{1-p}t_0^{(p-1)/2}.
\end{eqnarray}
This implies that 
\begin{eqnarray}
t_0 = \left(\frac{2b}{(3-p)\lambda\Vert W_p'\Vert^{1-p}}\right)^{2/(p-1)}.
\end{eqnarray}
By this and (4.1), we have 
\begin{eqnarray}
a =  \frac{p-1}{2}\lambda \Vert W_p'\Vert^{1-p}t_0^{(p-3)/2} =  \frac{p-1}{2}\lambda^{2/(p-1)}
\Vert W_p'\Vert^{-2}\left(\frac{2b}{3-p}\right)^{(p-3)/(p-1)}.
\end{eqnarray}
Therefore, we see that the equation $y(t) = g(t)$ has one solution $t_0 > 0$ if (4.5) holds. Then by (3.1) and (3.2), we know that 
$u_{\lambda}(x):= t_0^{1/2}\Vert W_p'\Vert^{-1}W_p(x)$ is a unique solution to 
(1.1). 
Equally, if 
\begin{eqnarray}
a <  \frac{p-1}{2}\lambda \Vert W_p'\Vert^{1-p}t_0^{(p-3)/2} = 
\frac{p-1}{2}\lambda^{2/(p-1)}
\Vert W_p'\Vert^{-2}\left(\frac{2b}{3-p}\right)^{(p-3)/(p-1)},
\end{eqnarray}
then $y(t) = g(t)$ has exactly two solutions $t_1, t_2$, and (1.1) has exactly two solutions, 
and if 
\begin{eqnarray}
a >  \frac{p-1}{2}\lambda \Vert W_p'\Vert^{1-p}t_0^{(p-3)/2} = 
\frac{p-1}{2}\lambda^{2/(p-1)}
\Vert W_p'\Vert^{-2}\left(\frac{2b}{3-p}\right)^{(p-3)/(p-1)},
\end{eqnarray}
then the equation $y(t) = g(t)$ has no solutions, and (1.1) has no solutions. Thus the proof is complete. 
\qed 

\vspace{0.2cm}

\noindent
{\bf Proof of Theorem 1.4.} 
In what follows, let $p = 2$. We use (3.1) and (3.2) here. 
We calculate $\eta =W_2(0)$ and $\Vert W_2'\Vert$.  By (1.6) and the same argument as that to obtain 
(2.7), for $-1 \le x \le 0$, we have 
\begin{eqnarray}
\frac12W_2'(x)^2 + \frac13W_2(x)^3 = \frac13\eta^3. 
\end{eqnarray}
This implies that for $-1 \le x \le 0$, 
\begin{eqnarray}
W_2'(x) = \sqrt{\frac{2}{3}}\sqrt{\eta^3 - W_2(0)^3}.
\end{eqnarray}
By this, we obtain 
\begin{eqnarray}
1 &=& \int_{-1}^0\sqrt{\frac{3}{2}}\frac{W_2'(x)}{\sqrt{\eta^3 - W_2(x)^3}}dx 
\\
&=&  \sqrt{\frac{3}{2}}\int_0^\eta\frac{1}{\sqrt{\eta^3-\theta^3}}d\theta
\nonumber
\\
&=&  \sqrt{\frac{3}{2}}\eta^{-1/2}\int_0^1 \frac{1}{\sqrt{1-s^3}}ds.
\nonumber
\end{eqnarray}
By this and (1.11), we have 
\begin{eqnarray}
\eta = W_2(0) = Q_1 = \frac32 A_2^2.
\end{eqnarray}
By (1.6) and (4.9), we have 
\begin{eqnarray}
\Vert W_2'\Vert^2 &=& 2\int_{-1}^0 \sqrt{\frac{2}{3}}\sqrt{\eta^3 - W_2(x)^3}W_2'(x)dx
\\
&=& 2\sqrt{\frac{2}{3}}\int_0^\eta \sqrt{\eta^3-\theta^3}d\theta 
= 2\sqrt{\frac{2}{3}}\eta^{5/2}\int_0^1\sqrt{1-s^3}ds
\nonumber
\\
&=& \frac{9}{2}A_2^5B_2.
\nonumber
\end{eqnarray}
This implies 
\begin{eqnarray}
\Vert W_2'\Vert = \frac{3}{\sqrt{2}}A_2^{5/2}B_2^{1/2}:= Q_2.
\end{eqnarray}
Now we prove Theorem 1.4 (i). Put $p = 2$ in Theorem1.3 (b). Then if 
\begin{eqnarray}
a = \frac{1}{4b}\lambda^2\Vert W_2'\Vert^{-2},
\end{eqnarray}
namely, if $\lambda = 2\sqrt{ab}\Vert W_2'\Vert$, then (1.1) has exactly one solution. By (3.1) and 
(4.3), we have 
\begin{eqnarray}
t_0^{1/2} = \sqrt{\frac{b}{a}}, \enskip Q_2^{-1} = 2\lambda^{-1}\sqrt{ab}.
\end{eqnarray}
By this and (3.1), 
\begin{eqnarray}
u_\lambda(x) = t_0^{1/2}Q_2^{-1}W_2(x) = 2b\lambda^{-1}W_2(x)
\end{eqnarray}
is the solution. By this and (4.11), we have 
\begin{eqnarray}
\eta &=& u_\lambda(0) = 2b\lambda^{-1}Q_1.
\end{eqnarray}
This implies Theorem 1.4 (i). We next prove Theorem 1.4 (ii). By Theorem1.3 (a), 
if $\lambda >2\sqrt{ab}Q_2$, then (1.1) has exactly two solutions. 
In this case, by (3.2), we have 
\begin{eqnarray}
t_1^{1/2}(\lambda) &:=& \frac{\lambda Q_2^{-1} + \sqrt{\lambda^{2}Q_2^{-2}-4ab}}{2a}, 
\\
t_2^{1/2}(\lambda) &:=& \frac{\lambda Q_2^{-1} - \sqrt{\lambda^{2}Q_2^{-2}-4ab}}{2a}. 
\end{eqnarray}
Then by (3.1), we have 
\begin{eqnarray}
u_{\lambda,1}(x) &=& \frac{\lambda Q_2^{-1} + \sqrt{\lambda^{2}Q_2^{-2}-4ab}}{2a}Q_2^{-1}W_2(x),
\\
u_{\lambda,2}(x) &=& \frac{\lambda Q_2^{-1} - \sqrt{\lambda^{2}Q_2^{-2}-4ab}}{2a}Q_2^{-1}W_2(x).
\end{eqnarray}
We put $x = 0$ in (4.20) and (4.21). Then by (4.11), we obtain (1.13) and (1.14). Thus the proof is complete. \qed

\section{Proof of Theorem 1.7}

{\bf Proof of  Theorem 1.7 (i).} Let $p = 3$. If we have $t_\lambda$ in (3.1) for some $\lambda > 0$, then we have a 
solution of (1.1) for $p = 3$. We put $p = 3$ in (3.1). Then we have 
\begin{eqnarray}
at + b = \lambda \Vert W_3'\Vert^{-2}t.
\end{eqnarray}
Namely, if $t_\lambda = \frac{b\Vert W_3'\Vert^2}{\lambda - a\Vert W_3'\Vert^2}$ exists, then 
we have a unique solution of (1.1) for $p = 3$. Assume that
\begin{eqnarray}
\lambda > a\Vert W_3'\Vert^2.
\end{eqnarray}
Then we have the unique solution 
\begin{eqnarray}
u_\lambda(x) =  \sqrt{\frac{b\Vert W_3'\Vert^2}{\lambda - a\Vert W_3'\Vert^2}}
\Vert W_3'\Vert^{-1}W_3(x).
\end{eqnarray}
Now we calculate $\eta = \Vert W_3'\Vert$. By (1.1) and the same argument as that to obtain (4.8), we have 
\begin{eqnarray}
\frac12W_3'(x)^2 + \frac14W_3(x)^4 = \frac14\eta^4. 
\end{eqnarray}
By this, for $-1 \le x \le 0$, we have 
\begin{eqnarray}
W_3'(x) = \frac{1}{\sqrt{2}}\sqrt{\eta^4-W_3(x)^4}.
\end{eqnarray}
By this, we have 
\begin{eqnarray}
\Vert W_3'\Vert^2 &=& 2\int_{-1}^0  \frac{1}{\sqrt{2}}\sqrt{\eta^4-W_3(x)^4}W_3'(x)dx 
\\
&=& \sqrt{2}\int_0^\eta\sqrt{\eta^4 - \theta^4}d\theta 
= \sqrt{2}\eta^3 \int_0^1 \sqrt{1-s^4}dx 
\nonumber
\\
&=& \sqrt{2}\eta^3B_3.
\nonumber
\end{eqnarray}
By (5.5), we have 
\begin{eqnarray}
1 &=& \int_{-1}^0 \sqrt{2}\frac{W_3'(x)}{\sqrt{\eta^4-W_3(x)^4}}dx
\\
&=& \sqrt{2}\int_0^\eta \frac{1}{\sqrt{\eta^4-\theta^4}}d\theta 
= \sqrt{2}\eta^{-1}\int_0^1 \frac{1}{\sqrt{1-s^4}}ds 
\nonumber
\\
&=& \sqrt{2}\eta^{-1}A_3.
\nonumber
\end{eqnarray}
By (5.6) and (5.7), we have 
\begin{eqnarray}
\Vert W_3'\Vert^2 = 4A_3^{3}B_3.
\end{eqnarray}
By this and (5.2), if $\lambda > 4aA_3^3B_3$, then we have the unique solution. We know 
from (1.17) that $\mu_1 = 4A_3^3B_3$. Then the argument in the proof of Theorem 1.2 is 
also available for the case $p = 3$, and we also obtain (1.7) for $p = 3$. This implies 
(1.20). Thus the proof of Theorem 1.7 (i) is complete. \qed 

\vspace{0.2cm}

\noindent
{\bf Proof of  Theorem 1.7 (ii).} Let $p = 1$. Then (1.1) is the linear eigenvalue problem 
with positive solution. Therefore, we have 
\begin{eqnarray}
\frac{\lambda}{b + a\Vert u'\Vert^2} = \lambda_1 = \frac{\pi^2}{4}. 
\end{eqnarray}
We look for the solution $u_\lambda(x) = C\varphi_1(x)$, where $C > 0$ is a constant and 
$\varphi_1(x):= \cos \frac{\pi}{2}x$. 
We substitute $C\varphi_1(x)$ into (5.9) to obtain 
\begin{eqnarray}
C
= \frac{4}{\pi^2\sqrt{a}}\sqrt{\lambda - \frac{\pi^2}{4}b}.
\end{eqnarray}
Therefore, if $\lambda > b\lambda_1$, then we have a unique solution 
$u_\lambda(x)$ of (1.1) such as
\begin{eqnarray}
u_\lambda(x) &=& 
\frac{4}{\pi^2\sqrt{a}}\sqrt{\lambda - \frac{\pi^2}{4}b}
\cos\frac{\pi}{2}x. 
\end{eqnarray}
Thus the proof of Theorem 1.7 (ii) is complete. \qed


\begin{thebibliography}{20}
\labelsep=1em\relax


\bibitem[1]{1} 
C.O. Alves, F.J.S.A. Corr\^{e}a, T.F. Ma, Positive solutions for a quasilinear elliptic equation of Kirchhoff type. Comput. Math. Appl. 49 (1), (2005), 85--93.
\bibitem[2]{2} {A. Ambrosetti, H. Brezis, G. Cerami}, 
{Combined effects of concave and convex 
nonlinearities in some elliptic problems}, 
J. Funct. Anal. {122} (2) (1994), 519--543. 
\bibitem[3]{3} X. Cao, G. Dai, Spectrum, global bifurcation and nodal solutions to Kirchhoff-type equations, Electron. J. Differential Equations 2018, Paper No. 179, 10 pp.
\bibitem[4]{4} B. Cheng, New existence and multiplicity of nontrivial solutions for nonlocal elliptic Kirchhoff type problems. J. Math. Anal. Appl. 394 (2) (2012),  488--495.
\bibitem[5]{5} D. Daners, J. L\'opez-G\'omez, Global dynamics of generalized logistic equations. Adv. Nonlinear Stud. 18 (2) (2018), 217--236. 
\bibitem[6]{6} B. Gidas, W. M. Ni, L. Nirenberg, Symmetry and related properties via the maximum principle. Comm. Math. Phys. 68 (3) (1979), 209--243.
\bibitem[7]{7} {T. Laetsch}, {
The number of solutions of a nonlinear 
two point boundary value problem}, 
Indiana Univ. Math. J. {20} (1) 1970, 1--13.
\bibitem[8]{8} 
Z. Liang, F. Li, J. Shi, Positive solutions of Kirchhoff-type non-local elliptic equation: a bifurcation approach. Proc. Roy. Soc. Edinburgh Sect. A 147 (4) (2017), 875--894.
\bibitem[9]{9}  K. Perera, Z. Zhang, Nontrivial solutions of Kirchhoff-type problems via the Yang index. J. Differential Equations 221 (1) (2006), 246--255. 
Springer, New York, 2002. 
\bibitem[10]{10} M. Sun, Z. Yang, H. Cai, Nonexistence and existence of positive solutions for the Kirchhoff type equation. Appl. Math. Lett. 96 (10) (2019), 202--207. 
\bibitem[11]{11} 
W. Wang, W. Tang, Bifurcation of positive solutions for a nonlocal problem. Mediterr. J. Math. 13 (6) (2016), 3955--3964.


\end{thebibliography}
\end{document}